\documentclass[11pt,a4paper]{article}
\usepackage{graphicx}
\usepackage{amsfonts}
\usepackage{latexsym}
\usepackage{epsfig}

\usepackage{mathpple}

\newcommand{\beql}[1]{\begin{equation}\label{#1}}
\newcommand{\eeq}{\end{equation}}
\newcommand{\comment}[1]{}

\newcommand{\eqref}[1]{{\rm (\ref{#1})}}

\newcommand{\Qed}{\ \\\mbox{$\Box$}}

\newcommand{\Set}[1]{{\left\{{#1}\right\}}}

\newcommand{\ZZ}{{\mathbb Z}}

\newcounter{open}
\newcounter{dfn}
\def\thedfn{\arabic{dfn}}

\newcounter{obs}
\def\theobs{\arabic{obs}}

\newcounter{thm}
\newcounter{othm}
\def\theothm{\Alph{othm}}
\newenvironment{othm}{
  \sf
  \vskip 0.10in
  \refstepcounter{othm}
  \noindent{\bf Theorem\ \theothm}
}

\newcounter{mysec}
\newcounter{mysubsec}[mysec]
\newtheorem{theorem}{Theorem}
\newtheorem{corollary}{Corollary}
\newtheorem{lemma}{Lemma}

\newtheorem{definition}{Definition}

\begin{document}

\title{Translational tilings of the integers with long periods}
\author{Mihail N. Kolountzakis}
\date{October 2002}
\maketitle

\begin{abstract}
Suppose that 
$A \subseteq \ZZ$, is a finite set of integers of diameter $D=\max A - \min A$.
Suppose also that $B \subseteq \ZZ$
is such that $A\oplus B = \ZZ$, that is each $n\in\ZZ$ is uniquely expressible as
$a+b$, $a\in A$, $b\in B$.
We say then that $A$ tiles the integers if translated at the locations $B$
and it is well known that $B$ must be a periodic set in this case.
Here we study the relationship between the diameter of $A$
and the least period $T$ of $B$.
We show that $T \le C \exp(C \sqrt D \log D \sqrt{\log\log D})$
and that we can have
$T \ge C D^2$, where $C>0$ denotes absolute constants.
\end{abstract}

\section{Introduction and results}

Let us start by defining translational tiling in any abelian group.
\begin{definition}\label{def:tiling}
Suppose $G$ is an abelian group, and $A,B \subseteq G$.
We say that $A$ tiles the set $E \subseteq G$ with $B$, and write $A\oplus B=E$,
if every element $g\in E$ admits a unique representation
$g=a+b$, $a\in A, b\in B$.
\end{definition}
It has long been known (see e.g. \cite{newman})
that, when $G=\ZZ$ is the integer group, $A \oplus B = \ZZ$
and $A \subseteq \ZZ$ is a finite set, then $B$ is necessarily periodic.
\begin{definition}
Suppose $B$ is a subset of an abelian group $G$.
We call $B$ periodic if there exists $t \in G\setminus\Set{0}$
such that $B = B+t$.
If $G = \ZZ$ we denote by ${\cal P}(B)$ the least such period $t>0$.
\end{definition}
We are interested in the case when we know the diameter $D$ of $A$,
which we assume to be able to tile $\ZZ$, and ask how large can the period of
a tiling by $A$ be.
The main quantity of interest is defined below.
\begin{definition}\label{def:quantity}
If $D$ is a positive integer we write ${\cal T}(D)$
for the largest integer $k$ such that there exists
$A \subseteq \Set{0,\ldots,D}$ and $B \subseteq \ZZ$ with
${\cal P}(B) \ge k$ and $A \oplus B = \ZZ$.
\end{definition}
It was proved in \cite{newman} that whenever $A$ is finite
and $A\oplus B = \ZZ$ then ${\cal P}(B) \le 2^D$, where
$D$ is the diameter of the set $A$.
It has also been noted by some authors (e.g. \cite{coven-meyerowitz})
that no tilings $A \oplus B = \ZZ$ are known with
${\cal P}(B) > 2D$ and $A \subseteq \Set{0,\ldots,D}$.
So apparently the state of knowledge regarding the function ${\cal T}(D)$ is
\beql{state}
2D \le {\cal T}(D) \le 2^D.
\eeq 

In this paper we improve these bounds somewhat.
\begin{theorem}\label{th:main}
There are absolute constants $c_1, c_2, c_3 > 0$ such that, for $D\ge 1$,
\beql{main}
c_1 D^2 \le {\cal T}(D) \le c_2 \exp(c_3 \sqrt D \log D \sqrt{\log\log D}).
\eeq
\end{theorem}
Theorem \ref{th:main} follows directly from Theorems \ref{th:lower-bound} and
\ref{th:upper-bound} below.
\begin{theorem}\label{th:lower-bound}
For $D\ge 1$ 
there are tilings $X\oplus Y = \ZZ$ with
$X\subseteq\Set{0,\ldots,D}$ and
$$
{\cal P}(Y) \ge C D^2,
$$
where $C > 0$ is a constant.
\end{theorem}
\begin{theorem}\label{th:upper-bound}
If $A\oplus B = \ZZ$ is a tiling and $A \subseteq \Set{0,\ldots,D}$ then
$$
{\cal P}(B) \le C \exp(C \sqrt D \log D \sqrt{\log\log D}),
$$
where $C$ denotes absolute positive constants.
\end{theorem}
{\em Remarks. }
\begin{enumerate}
\item
We have not been able to construct a set $X$ which tiles
{\em only} with long periods. This would be very interesting as this
could be a finitary and translational analog of {\em aperiodic tiling},
the ability that is of a {\em set} of tiles to tile space but only
aperiodically \cite{grunbaum-shephard}.
\item
There is a theorem of de Bruijn \cite[Theorem 1]{debruijn} which can be used to prove
Theorem \ref{th:lower-bound}.
We think our proof is more intuitive
(but a more general problem is studied in \cite{debruijn}).
\item
The largest part of the proof of Theorem \ref{th:upper-bound} appears
also in an argument by Granville in \cite{konyagin-laba}.
\end{enumerate}

\noindent
{\bf Notation: } We write $[n] = \Set{0,\ldots,n-1}$.
The same set considered as an additive group mod $n$ is denoted by $\ZZ_n$.

Tiling the integers with a long period is equivalent to tiling a long finite
cycle in a non-periodic way, as the following easy lemma claims.
\begin{lemma}\label{lm:prelim}
Suppose $A \subseteq [M]$, $M$ a positive integer, and that $B \subseteq \ZZ$
is $M$-periodic: $B = B' + M\ZZ$, where $B' \subseteq [M]$.
Then $A\oplus B = \ZZ$ and $M = {\cal P}(B)$ if and only if
$A \oplus B' = \ZZ_M$ and $B'=B \bmod M$ is not periodic mod $M$.
\end{lemma}
{\bf Proof.}
Suppose that $A \oplus B = \ZZ$. It follows that
$A \oplus (B' + M\ZZ) = \ZZ$ hence $A \oplus B' = \ZZ/m\ZZ = \ZZ_M$.
If $t \in \ZZ_M\setminus\Set{0}$ and $B' = B'+t$ then
$B = (B'+t)+M\ZZ = B+t$ which contradicts the fact that $M$ is a minimal
period for $B$.
Conversely, suppose that $A\oplus B' = \ZZ_M$ and that $B'$ is not periodic mod $M$.
Then $A\oplus B = A \oplus B' + M\ZZ = [M]\oplus\ M\ZZ = \ZZ$ is a tiling of
period $M$. Furthermore, by the previous argument if $B$ has a smaller period $t \in [M]$
it follows that $B'$ has $t \bmod M$ as a period, a contradiction.
\Qed

However, it is not possible to have non-repetitive tilings of {\em intervals}
with tiles which are very short compared to the interval.
Essentially, in any tiling of an interval
by a tile of diameter $D$ the length of the interval is at most $2D$,
otherwise there is a smaller sub-tiling.
This fact is probably responsible for many researchers being
unaware of complicated tilings (such as those claimed by Theorem \ref{th:lower-bound}),
as one tends to think of tilings of the integers which are created
by first tiling an interval and then repeating that interval.
\begin{theorem}\label{th:interval}
If $0\in A\cap B$, $A \oplus  B = [n]$, $n>1$,
and $\max A > \max B$ then $A$ is periodic mod $n$.
\end{theorem}
{\em Remark.} Under the assumptions of Theorem \ref{th:interval} $\max A \neq \max B$,
so one of them is larger.
\begin{corollary}\label{cor:interval}
If $A \subseteq \Set{0,\ldots,D}$, $A\oplus B = [n]$ is a tiling and $n>2D$
there is $t<n$ such that $A \oplus (B\cap [t])$ is already a tiling of an interval.
\end{corollary}
{\bf Proof of Corollary \ref{cor:interval}.}
Since $n = \max A + \max B$ it follows that $\max B > \max A$
and, by Theorem \ref{th:interval}, $B = B + t \bmod n$, for some $t\in\Set{1,\ldots,n-1}$,
$t \mid n$.
From this it follows easily that $A \oplus  (B \cap [t]) = [t]$ is also a tiling.
\Qed

\section{Proofs}
In what follows $C$ denotes an absolute positive constant, not necessarily
the same in all its occurences.

\noindent
{\bf Proof of Theorem \ref{th:lower-bound}.}
By Lemma \ref{lm:prelim} it suffices to construct a set $A\subseteq\Set{0,\ldots,D}$
and a non-periodic set $B \subseteq \ZZ_M$, such that $A\oplus B = \ZZ_M$
is a tiling and $M \ge C D^2$.

Assume that $D$ is large and
pick $p$ and $q$ two primes in the range $D/2 \le p, q \le 2 D$, for example.
Let $M = 2 \cdot 3 \cdot 5 \cdot p \cdot q$, so that $M \sim D^2$.
The group $\ZZ_M$ is isomorphic to $\ZZ_{3p} \times \ZZ_{5q} \times \ZZ_2$, and
we visualize it as a parallelogram $Q$ in three dimensions with $\ZZ_{3p}$ along the $x$-axis,
$\ZZ_{5q}$ along the $y$-axis and $\ZZ_2$ along the $z$-axis.

\begin{figure}[h]
\begin{center}\input{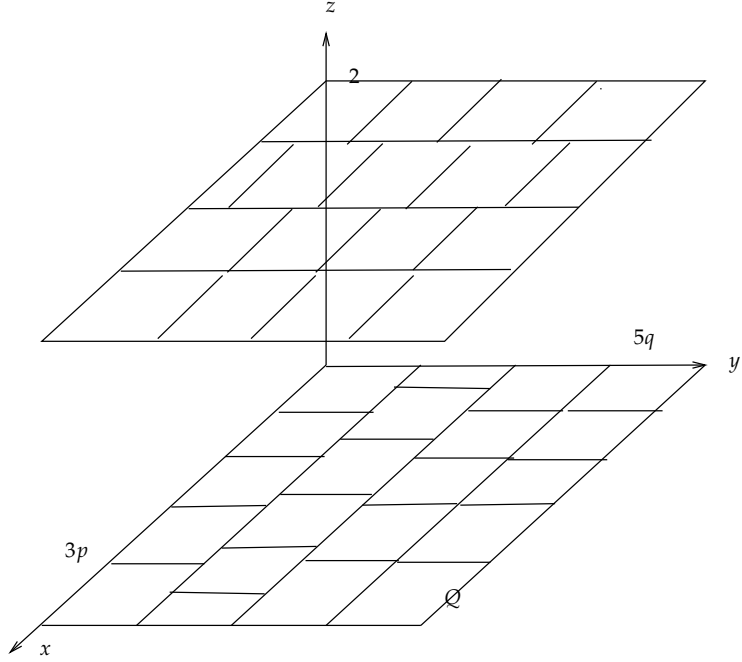}\end{center}
\caption{A non-periodic tiling by $A$}
\label{fig:tiling}
\end{figure}

We define $A$ as being the $3\times 5$ rectangle on the $xy$-plane with corner at $0$,
as is shown in Figure \ref{fig:cube}.
\begin{figure}[h]
\begin{center}\input{cube.pstex_t}\end{center}
\caption{The group $\ZZ_{3p} \times \ZZ_{5q} \times \ZZ_2$}
\label{fig:cube}
\end{figure}
Now we observe that $A$ can tile the group $Q$ non-periodically,
as shown in Figure \ref{fig:tiling}.
To obtain the tiling shown in Figure \ref{fig:tiling}, first tile $Q$ by $A$ in
the usual way, i.e.\ with $B$ equal to the subgroup generated by the elements
of $Q$: $(3,0,0)$, $(0,5,0)$ and $(0,0,1)$. Then, in the bottom level of the tiling
shift one of the columns along the $x$ direction, say the second one as shown,
by the vector $(1,0,0)$, and shift the second row in the upper level (one which
is parallel to the $y$ direction) by the vector $(0,1,0)$.

The set $B$ is therefore given by $B = L \cup U$ (lower and upper part) and
$L = \Set{(a,b,0)}$, where $b$ takes all values in $<5> \subseteq \ZZ_{5q}$ and
$a$ takes all values in $<3> \subseteq \ZZ_{3p}$, except when $b=5$ in which case
$a$ takes all values in $<3>+1 \subseteq \ZZ_{3p}$.
Similarly $U = \Set{(a,b,1)}$, where $a$ takes all values in $<3> \subseteq \ZZ_{3p}$
and $b$ takes all values in $<5> \subseteq \ZZ_{5q}$, except when $a=3$, in which
case $b$ takes all values in $<5>+1 \subseteq \ZZ_{5q}$.

It is clear that $B$ is not periodic. Indeed, if $t=(i,j,k) \in Q\setminus\Set{0}$ is
a period then it must belong to $B$ as $0 \in B$. If $t \in L$ then it must belong to the
group $<(3,0,0)>$ but $U$ shifted by any of these non-zero elements does not go into $U$,
because the second row shifted by such an element does not go into another row.
If $t \in U$ then $U+t$ ought to be equal to $L$ but it is clearly not.

Fix now a group isomorphism $\psi: Q \to \ZZ_M$ and define
$$
X = \psi(A),\ \ Y = \psi(B).
$$
Clearly $X\oplus Y=\ZZ_M$ is a tiling and $Y$ is not periodic,
as both these properties are invariant under group isomorphisms.
It remains to find the diameter of $X$, which of course depends on $\psi$.
Choose then the isomorphism
$$
\psi(i,j,k) = i(2\cdot5q)+j(2\cdot3p)+k(3p5q) \bmod M,\ \ (i=0,\ldots,3p-1, j=0,\ldots,5q-1, k=0,1),
$$
which maps $A$ into the set $X=\Set{0,\ldots,2\cdot5q+4\cdot3p}$.
This set $X$ has diameter $\le C D$ whereas $M \ge C D^2$.
\Qed

\noindent
{\em Remark.} In our example the tile $X$ has constant size $15$.

\noindent
{\bf Proof of Theorem \ref{th:upper-bound}.}
Suppose that $A \subseteq \Set{0,\ldots,D}$, $D<M$, and that $A\oplus B = \ZZ_M$,
for some $B \subseteq \ZZ_M$. We will show that if $M$ is large then $B$ is periodic.

For $B$ to be periodic we must show the existence of $t \in \ZZ_M\setminus\Set{0}$
such that $B = B+t$. Equivalently, writing for any finite set of integers $E$
$$
E(x) = \sum_{n \in  E} x^n,
$$
we must ensure that
$$
B(x) = x^t B(x) \bmod x^M-1.
$$
This is equivalent to
\beql{periodic-condition}
x^M-1 \mid (x^t-1) B(x).
\eeq

On the other hand, the fact that $A\oplus B$ is a tiling of $\ZZ_M$ is equivalent to
$$
A(x)B(x) = 1+x+x^2+\cdots+x^{M-1} \bmod x^M-1,
$$
which is equivalent to
\beql{tiling-condition}
x^M-1 \mid A(x)B(x) - {x^M-1 \over x-1},
\eeq
which implies that all $M$-th roots of unity except $1$
are roots of $A(x)$ or $B(x)$.
Equivalently, we must have that for each $d \mid M$, $d>1$, the
cyclotomic polynomial $\Phi_d(x)$ divides $A(x)$ or $B(x)$. 
And, similarly, \eqref{periodic-condition} is equivalent to
having all cyclotomic polynomials $\Phi_d(x)$, for $d \mid M$, $d>1$,
divide $x^t-1$ or $B(x)$.

Hence, for $t$ to be a period of $B$, it suffices that all cyclotomic
$\Phi_d(x)$, for $d \mid M$, $d>1$, that divide $A(x)$ also divide $x^t-1$.
Let now $\Phi_{s_1}(x),\ldots,\Phi_{s_k}(x)$ be all cyclotomic polynomials
$\Phi_s(x)$ with $s>1$ that divide $A(x)$, written {\em once} each
and numbered so that $1<s_1<s_2<\cdots<s_k$.
Since $\deg\Phi_s = \phi(s)$ (the Euler function) it follows that
\beql{degree-inequality}
\phi(s_1)+\cdots+\phi(s_k) \le \deg{A(x)} \le D.
\eeq
But \cite[p.\ 267]{hardy-wright}
$$
\phi(n) \ge C {n \over \log\log n}
$$
which implies $\log n \le C \log\phi(n)$, hence
\beql{s-sum}
s_1+\cdots+s_k \le C D \log\log s_k \le C D \log\log\phi(s_k) \le C D \log\log D.
\eeq
From the fact that the $s_i$ are
{\em different integers} and \eqref{s-sum} it
follows that $k \le C \sqrt{D \log\log D}$
hence
\beql{prod-bound}
\prod_{i=1}^k s_i \le (D \log\log D)^k \le \exp(C \sqrt D \log D \sqrt{\log\log D}).
\eeq
Define now $t = \prod_{i=1}^k s_i$,
so that all cyclotomic polynomials that divide $A(x)$ are also
divisors of $x^t-1$.
 From \eqref{prod-bound}, it follows that \eqref{periodic-condition}
holds, and therefore $t$ is a period of $B$, if only
$t<M$.
But $B$ is assumed to be non-periodic, so
$$
M \le C \exp(C \sqrt D \log D \sqrt{\log\log D}),
$$
as we had to prove.
\Qed

\noindent
{\bf Proof of Theorem \ref{th:interval}.}
We use the following result by Long \cite[Lemma 1]{long}.
\begin{othm}\label{th:long}
Suppose $n>1$ and $A\oplus B = [n]$.
Then there is an integer $m\ge 2$, $m\mid n$, and sets $E, D \subseteq [n/m]$, such that
$E\oplus D = [n/m]$, $0 \in E \cap D$, and
$$
\Set{A, B} = \Set{mE\oplus [m], mD}.
$$
\end{othm}

Suppose now that $n>1$ is minimal such that Theorem \ref{th:interval} fails
for $n$ and let $A, B \subseteq [n]$, $0\in A \cap B$, be such that $\max A > \max B$,
$A\oplus B = [n]$ and $A$ is not periodic mod $n$.
Now Theorem \ref{th:long} applies and let $m, E, D$ be as in the conclusion of it.
If $m=n$ it follows that $A = [n], B = \Set{0}$, so that $A$ is periodic mod $n$.
If $m<n$ then $n>n/m>1$ hence Theorem \ref{th:interval} holds for the integer $n/m$.
There are now two cases: (i) $\max E > \max D$, (ii) $\max D > \max E$.

In case (i) $\max(mE+[m]) > \max(mD)$ hence $A = mE+[m]$.
Since Theorem \ref{th:interval} holds for the integer $n/m$ it follows that
$E$ is periodic mod $n/m$ with some period $t$, which implies that $A$ is periodic
mod $n$ with period $mt$.

In case (ii) $\max(mD) > \max(mE+[m])$ hence $A = mD$.
By Theorem \ref{th:interval} again we get that $D$ is periodic mod $n/m$ with some period
$t$, and $A$ is therefore periodic mod $n$ with period $mt$.

In both cases we conclude that our theorem applies for $n$, hence it's true for all $n$.
\Qed

\noindent
\ \\
{\bf Bibliography}
\\

\noindent
{\sc\small
Department of Mathematics, University of Crete, Knossos Ave.,
714 09 Iraklio, Greece.\\
E-mail: {\tt mk@fourier.math.uoc.gr}
}

\end{document}